\documentclass[12pt]{amsart}

\usepackage[table,xcdraw]{xcolor}
\usepackage{amssymb}
\usepackage{amsmath}
\usepackage{amsfonts}
\usepackage{amsxtra}
\usepackage{amsthm}
\usepackage{tikz}
\usepackage{tikz-cd}
\usepackage{float}
\usepackage{booktabs}
\usepackage{color}
\usepackage{xcolor}
\usepackage{graphicx}
\usepackage{enumitem}
\usepackage{venndiagram}
\usetikzlibrary{positioning,shapes,fit,arrows}
\usepackage{multirow}
\usepackage{multicol}
\usepackage{hyperref}

\newtheorem{theorem}{Theorem}[section]
\newtheorem{proposition}[theorem]{Proposition}
\newtheorem{corollary}[theorem]{Corollary}
\newtheorem{lemma}[theorem]{Lemma}
\theoremstyle{definition}
\newtheorem{exmp}[theorem]{Example}
\newtheorem{definition}[theorem]{Definition}
\newtheorem{remark}[theorem]{Remark}

\usepackage{xcolor}
\definecolor{redue}{RGB}{147, 6, 3}

\newcommand{\abf}{\mathbf{a}}

\newcommand{\Id}{\mathrm{Id}}

\newcommand{\Md}{\mathrm{Med}}
\newcommand{\HomM}{\mathrm{Hom}(M,N)}

\newcommand{\Open}{\mathrm{Open}}

\def\z{\mathbb{Z}}
\def\q{\mathbb{Q}}
\def\r{\mathbb{R}}

\def\inv{^{-1}}
\def\fc{\mathcal{F}}
\newcommand{\Pd}{\mathrm{PreMed}}
\newcommand{\PS}{\mathrm{PreSh}(X)}
\newcommand{\supp}{\mathrm{supp}}
\begin{document}
	
\title{Bridging Meadows and Sheaves}
\author{João Dias, Bruno Dinis, and Pedro Macias Marques}
\date{}

 \address[João Dias]{Departamento de Matemática, Universidade de Évora}
\email{joao.miguel.dias@uevora.pt}

\address[Bruno Dinis]{Departamento de Matemática, Universidade de Évora}
\email{bruno.dinis@uevora.pt}

\address[Pedro Marques]{Departamento de Matemática, Universidade de Évora}
\email{pmm@uevora.pt}



	\begin{abstract}
		We bridge sheaves of rings over a topological space with common meadows (algebraic structures where the
inverse for multiplication is a total operation). More specifically, we show that the subclass of pre-meadows with $\abf$, coming from the lattice of open sets of a topological space $X$, and presheaves over $X$ are the same structure. 
Furthermore, we provide a construction that, given a sheaf of rings $\mathcal{F}$ on $X$ produces a common meadow as a disjoint union of elements of the form $\mathcal{F}(U)$ indexed over the open subsets of $X$. We also establish a correspondence between the process of going from a presheaf to a sheaf (called sheafification) and the process of going from a pre-meadow with $\abf$ to a common meadow.
	\end{abstract}
 \maketitle
 
	\section{Introduction}

Meadows were originally introduced by Bergstra and Tucker in \cite{Bergstra2006}, as abstract data types, given by equational axiomatizations, in order to obtain simple and easy to automate in formal reasoning, term rewriting systems \cite{Bergstra2008,BP(20),10.1145/1219092.1219095,bergstra2020arithmetical,bergstra2020arithmetical,bergstra2023axioms}. 

Meadows can also be seen as algebraic structures with two operations, addition and multiplication, where both operations have inverses which are total. From this point of view, the main novelty of these structures is the possibility to divide by zero. This turns out to be a fruitful line of research, which has seen some development recently. See for instance  \cite{Dinis_Bottazzi,Dias_Dinis(23),Dias_Dinis(24),Dias_Dinis(Art),DD(Flasque),Bergstra_Tucker_2024, bergstra2024ringscommondivisioncommon}. 

There are several different classes of meadows but for our purposes two of them are particularly relevant: \emph{common meadows}, introduced by Bergstra and Ponse in \cite{Bergstra2015} and \emph{pre-meadows with $\abf$}, introduced by the first two authors in \cite{Dias_Dinis(23)}.  Pre-meadows with $\abf$ differ from common meadows as the former do not require the existence of an inverse function. The two structures are nevertheless related, since in both of them there is an element, denoted $\abf$ (sometimes also denoted $\bot$), which is the inverse of zero, and is otherwise absorbent for both operations.  Moreover, if $P$ is a pre-meadow with $\abf$ such that $0\cdot P$ is finite and the partial order in $0\cdot P$ is a total order, then it is in fact a common meadow. 

A result from \cite{Dias_Dinis(23)} (see Theorem~\ref{T:DirectedLattice} below)  shows that pre-meadows with $\abf$ (and common meadows) can be decomposed as disjoint unions of rings. This fact allows one to adapt many results from ring theory into the setting of meadows. The paper \cite{Dias_Dinis(23)} also contains a section in which meadows are considered from a categorical point of view. In this paper, we take this view a step forward, and bridge meadows with algebraic topology, by establishing a connection between meadows and sheaves. 

We show that a subclass of pre-meadows with $\abf$, coming from the lattice of open sets of a topological space $X$ and presheaves over $X$ are the same structure. 
Furthermore, given a sheaf of rings $\mathcal{F}$ on $X$ one is able to construct a common meadow as a disjoint union of elements of the form $\mathcal{F}(U)$ indexed over the open subsets of $X$. We also establish a correspondence between the process of going from a presheaf to a sheaf (called sheafification) and the process of going from a pre-meadow with $\abf$ to a common meadow.

 The paper is structured as follows. We start by recalling some preliminary notions and results related to sheaves (Section~\ref{S:Sheaves}) and common meadows (Section~\ref{S:Meadows}). In Section~\ref{S:Bridge}, we show results that establish a connection between sheaves and common meadows.  This connection is established by showing first, in Section~\ref{S:SheavestoMeadows}, how to construct a common meadow whenever a sheaf of rings over a topological space is given. Then, in Section~\ref{S:Relation with Sheafification},  we show that the process of going from pre-meadows with $\abf$ to common meadows can be identified with the sheafification  of a presheaf. However, if one starts with a common meadow, this process may give an essentially  different common meadow. So, in Section~\ref{S:Isomorphism}, we give conditions for the resulting common meadow to be isomorphic to the original one. Some concluding remarks are left for Section~\ref{S:Conclusion}. 
	\section{Preliminaries}

	\subsection{Sheaves}\label{S:Sheaves}
We recall the definitions of presheaf and sheaf of rings (see Chapter~II, Section~1 in \cite{Hartshorne(77)} for a nice introduction to sheaves). Unless otherwise stated, all rings we consider are commutative.
	\begin{definition}
Let $X$ be a topological space. A \emph{presheaf of (commutative) rings} $\mathcal{F}$ on $X$ is a family of (commutative) rings indexed by the open sets of $X$ together with a family of homomorphisms of rings, indexed by pairs of open sets of $X$, such that:
		\begin{enumerate}
			\item For every open subset $U$ of $X$, the element $\fc(U)$ is a (commutative) ring.
			\item If $V$ and $U$ are open subsets of $X$ such that $V\subseteq U$, then there exists a ring homomorphism $\rho_{V,U}:\fc(U)\rightarrow \fc(V)$ such that
		\begin{enumerate}
			\item $\fc(\emptyset)=0$, the zero ring.
			\item $\rho_{U,U}$ is the identity.
			\item If $W\subseteq V\subseteq U$ are open subsets of $X$, then $\rho_{W,U}=\rho_{W,V}\circ\rho_{V,U}$.
		\end{enumerate}
		\end{enumerate}
	\end{definition}
Observe that we are writing $\rho_{V,U}$ to denote the restriction map from $\fc(U)$ to $\fc(V)$, where, for example, in \cite{Hartshorne(77)} the corresponding notation is $\rho_{U,V}$. The reason for this choice is that we are making this consistent with the usual notation for transition maps of common meadows.
	
	\begin{definition}
		If $X$ is a topological space and $\fc$ is a presheaf over $X$, then we say that $\fc$ is a \emph{sheaf} if for any open set $U$ and any open cover $\{U_i\}_{i\in I}$ of $U$ (i.e., a family of open sets such that  ${U=\bigcup_{i\in I}U_i}$), we have		
		\begin{enumerate}
			\item (\textit{Locality}) 
For any $s,t\in \fc(U)$, if $\rho_{U_i,U}(s)=\rho_{U_i,U}(t)$ for all $i\in I$ then $s=t$.
			\item (\textit{Glueing}) 
For any collection $\{s_i\}_{i\in I}$ satisfying ${s_i\in \fc(U_i)}$ for each ${i \in I}$, if we have $\rho_{U_i\cap U_j,U_i}(s_i)=\rho_{U_i\cap U_j,U_j}(s_j)$ for all ${i,j\in I}$ then there exists $s\in U$ such that $\rho_{U_i,U}(s)=s_i$, for all ${i \in I}$.
		\end{enumerate}
	\end{definition}
 \begin{exmp}\label{E:Not_a_Sheaf}
     Let $\underline{\z}$ be the constant presheaf  on the set of real numbers $\r$, where to each non-empty subset we associate the ring $\z$, and to the empty set the trivial ring. It is well-known that this presheaf is not a sheaf. Indeed, the glueing property fails, since if we take two disjoint open subsets, and different sections, they must coincide in the intersection however they do not come from the same section in the union.
 \end{exmp}
 
	The process of transforming a presheaf into a sheaf is known as \emph{sheafification}. In order to recall this notion we need the definition of stalk.
	
	\begin{definition}
Let $X$ be a topological space, let $\fc$ be a presheaf on $X$, and let $x\in X$. The \emph{stalk} of $\fc$ at $x$ is defined as:
		$$\fc_x=\{(s,U)\mid x\in U\text{ and } s\in \fc(U) \}/\sim,$$
		where $(s,U)\sim (t,V)$ if $\rho_{U\cap V,U}(s)=\rho_{U\cap V,V}(t)$. If ${x \in U}$ and ${s \in \fc (U)}$, we denote by $s_x$ the class of $(s,U)$ in $\fc_x$.
	\end{definition}
	
\begin{definition}
		Let $X$ be a topological space, and $\fc$ a presheaf on $X$.
        We say that an element ${(s_x)_{x \in U}}$ in the product $\prod_{x\in U}\fc_x$ satisfies condition $(\ast)$ if for any  ${x \in U}$, there exist an open set $V$, with ${x\in V\subseteq U}$,
		and an element ${t\in \fc(V)}$ such that for all ${y\in V}$ we have  ${t_y=s_y}$.
        We can then define the sheaf
		$$\fc^+(U)=\biggl\{\{s_x\}_{x\in U}\in\prod_{x\in U}\fc_x\mid \text{ condition } (\ast) \text{ holds} \biggr\}.$$
\end{definition}
   
	\begin{definition}
		Let ${f:X\to Y}$ be a continuous map. Then \emph{the direct image functor} $f_*:\mathrm{PreSheaves}(X)\to \mathrm{PreSheaves}(Y)$ is defined by
		$$(f_*(\fc))(U)=\fc(f\inv(U)),$$
		where $\fc$ is a presheaf in $X$ and $U$ is an open set in $Y$. 
	\end{definition}
	
\begin{remark}
\label{R:directimagecommutes}
We can easily check that given a continuous map ${f: X\to Y}$, the direct image functor $f_*$ takes sheaves to sheaves. Furthermore, in the case where $f$ is an isomorphism, it  commutes with sheafification. 
\end{remark}

	\subsection{Common Meadows}\label{S:Meadows}
We recall the definition of pre-meadow, pre-meadow with $\abf$, and common meadow, and some of their properties. We refer to \cite{Dias_Dinis(23)} for more details.
\begin{definition}
A \emph{pre-meadow} is a set $P$, having two distinct elements denoted $0_P$ and $1_P$, or $0$ and $1$ when $P$ is understood, and endowed with binary operations $+$ and $\cdot$, and a unary operation $-$ satisfying the following properties
    \begin{multicols}{2}
\begin{enumerate}
\item[$(P_1)$] $(x+y)+z=x+(y+z) $
\item[$(P_2)$] $x+y=y+x $ 
\item[$(P_3)$]  $x+0=x$ 
\item[$(P_4)$] $x+ (-x)=0 \cdot x$
\item[$(P_5)$] $(x \cdot y) \cdot z=x \cdot (y \cdot z)$ 
\item[$(P_6)$]  $x \cdot y=y \cdot x $
\item[$(P_7)$] $1 \cdot x=x$
\item[$(P_8)$] $x \cdot (y+z)= x \cdot y + x \cdot z$
\item[$(P_9)$] $-(-x)=x$
\item[$(P_{10})$]$0\cdot (x+y)=0\cdot x \cdot y$ \\
\end{enumerate}
\end{multicols}
\end{definition}

	\begin{definition}\label{ringsinameadow}
		Let $P$ be a pre-meadow and let ${z\in0\cdot P}$. Then $P_z$ denotes the set
\[
P_z=\{x\in P \mid 0\cdot x=z\}.
\]
	\end{definition}
	
	\begin{definition}
We say that $P$ is a \emph{pre-meadow with $\abf$} if there exists a unique element in $0\cdot P$ (denoted by $\abf$) such that $|P_\abf|=1$  and $x+\abf=\abf$, for all $x\in P$.
    \end{definition}
    
    \begin{definition}
    A \emph{common meadow} is a pre-meadow with $\abf$ equipped with an inverse function $(\cdot)\inv$ satisfying
    \begin{multicols}{2}
\begin{enumerate}
\item[$(M_1)$] $x \cdot x^{-1}=1 + 0 \cdot x^{-1}$
\item[$(M_2)$]$(x \cdot y)^{-1} = x^{-1} \cdot y^{-1}$
\item[$(M_3)$] $(1 + 0 \cdot x)^{-1} = 1 + 0 \cdot x $
\item[$(M_4)$] $ 0^{-1}=\abf$
\end{enumerate}
\end{multicols}
\end{definition}

Each pre-meadow has an associated directed lattice whose vertices are rings. The ring homomorphisms of the associated lattice are called \emph{transition maps}. 

Homomorphisms of meadows are defined as follows.
\begin{definition}\label{D:Morphism}
    Let $f:M\rightarrow N$ be a function. We say that $f$ is an \emph{homomorphism of (common) meadows} if $M,N$ are common meadows and for all $x,y\in M$
    \begin{enumerate}
        \item $f(x+y)=f(x)+f(y)$.
        \item $f(x\cdot y)=f(x) \cdot f(y)$.
        \item $f(1_M)=1_{N}$.        
    \end{enumerate}
\end{definition}

	\begin{definition}\label{D:LatticeRings}
		A \emph{directed lattice} of rings $\Gamma$ over a countable lattice $L$ consists of a family of commutative rings $\Gamma_i$ indexed by $i\in L$, such that $\Gamma_i$ is a unital commutative ring for all $i\in L\setminus \min(L)$ and $\Gamma_{\min(L)}$ is the zero ring, together with a family of ring homomorphisms $f_{j,i}:\Gamma_i\rightarrow\Gamma_j$ whenever $i>j$ such that $f_{k,j}\circ f_{j,i}=f_{k,i}$ for all $i>j>k$.
	\end{definition}
	
	We denote by $R^{\times}$ the set of invertible elements of the ring $R$. In \cite{Dias_Dinis(23)} the following connection between rings and meadows was shown.

\begin{theorem}\label{T:DirectedLattice}
        Let $\Gamma$ a directed lattice of rings over a lattice  $L$. 
          Then $M=\bigsqcup_{i\in L}\Gamma_i$, endowed with the operations 
          \begin{itemize}
                \item $x+_My=f_{i\wedge j,i}(x)+_{i\wedge j}f_{i\wedge j,j}(y)$,
                \item $x \cdot_M y=f_{i\wedge j,i}(x)\cdot_{i\wedge j}f_{i\wedge j,j}(y)$,
            \end{itemize}
            where $+_{i\wedge j}$ and $\cdot_{i\wedge j}$ are the sum and the product in $\Gamma_{i\wedge j}$, respectively,
             is a pre-meadow with $\abf$ and the lattice $0\cdot M$ is equivalent to $L$.\\ Moreover, $M$ is a common meadow if and only if, for all $i\in L$ and all $x\in \Gamma_i$ the set
        $$J_x=\{j\in L\mid f_{j,i}(x)\in \Gamma_j^{\times}\}$$
        has a maximum.  
\end{theorem}

	\begin{corollary}\label{C:DirectedLattice}
		Let $L$ be a lattice and $\Gamma$ a directed lattice of rings over $L$. Then, there exists $M=\bigsqcup_{i\in I}\Gamma_i$, a pre-meadow  with $\abf$, endowed with the operations as in Theorem~\ref{T:DirectedLattice} and the lattice $0\cdot M$ is equivalent to $L$. 
	\end{corollary}

\begin{remark}
By \cite[Theorem~2.5]{Dias_Dinis(23)}, $M_z$ is a ring with zero $z \in 0 \cdot M$ and identity ${1+z}$.
\end{remark}
\begin{definition}
Let $M$ be a pre-meadow. We define a partial order on $0\cdot M$ by 
\[
z \le z' \quad \text{if and only if} \quad z \cdot z' = z.
\] 
\end{definition}

	\begin{proposition}\label{P:Transitionmaps}
		Let $M$ be a meadow. If $z,z'\in 0\cdot M$ are such that $z\leq z'$, then the map
		\begin{align*}
			f_{z,z'}:M_{z'}&\rightarrow M_z\\
			x&\mapsto x+z
		\end{align*}
		is a ring homomorphism.
		
		Moreover, if $z,z',z''\in 0\cdot M$ are such that $z\leq z'\leq z''$, then ${f_{z,z'}\circ f_{z',z''}=f_{z,z''}}$.
	\end{proposition}

	\section{Meadows as Sheaves}\label{S:Bridge}

	In this section we show that, under certain conditions, pre-meadows and presheaves are essentially the same structure. In other words, there exist adjoint functors between the category of presheaves and the category of pre-meadows over open sets. Of particular interest are pre-meadows defined over the same lattice, that is, pre-meadows $M$ and $N$ such that the lattice $0\cdot M$ is isomorphic to $0\cdot N$.

 In Section~\ref{S:SheavestoMeadows} we show how to construct a common meadow whenever a sheaf of rings over a topological space is given. Then, in Section~\ref{S:Relation with Sheafification},  we show that the process of going from pre-meadows with $\abf$ to common meadows can be identified with the sheafification  of a presheaf. One may ask, if one already starts with a common meadow, whether this process gives the same common meadow. Alas, that is not the case. However, as shown in Section~\ref{S:Isomorphism}, it is possible to characterize when the resulting common meadow is  isomorphic to the original one.

 \subsection{From Sheaves to Meadows}\label{S:SheavestoMeadows}

	\begin{definition}
		The category  $\Pd$ of pre-meadows with $\abf$ is defined as follows:
		\begin{itemize}
			\item The \emph{objects} of $\Pd$ are pre-meadows with $\abf$.
			\item If $M$ and  $N$ are objects in $\Pd$ then the \emph{morphisms} from $M$ to $N$ are the elements of $\HomM$, the set of all homomorphisms of pre-meadows with $\abf$ from $M$ to $N$.
		\end{itemize}
	\end{definition}
	
	To each pre-meadow with $\abf$ we may associate a lattice $0\cdot M$. Conversely, given a lattice $L$ we can define the (non-empty) full subcategory $\Pd(L)$ of $\Pd$ consisting of the pre-meadows with $\abf$, say $M$, such that $0\cdot M$ is isomorphic to $L$. 
	
	\begin{exmp}
Let $L$ be the lattice
		
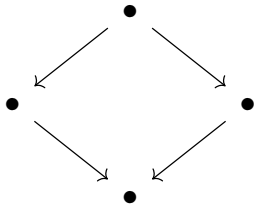
\[\begin{tikzcd}
			& \bullet \\
			\bullet && \bullet \\
			& \bullet
			\arrow[from=1-2, to=2-1]
			\arrow[from=1-2, to=2-3]
			\arrow[from=2-1, to=3-2]
			\arrow[from=2-3, to=3-2]
		\end{tikzcd}\]
Then the category $\Pd(L)$ consists of all meadows $M$ such that the associated directed lattice is:
		
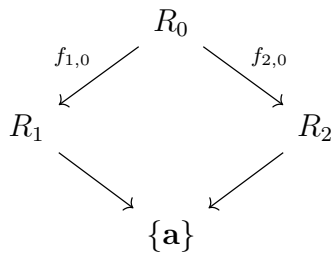
\[\begin{tikzcd}
			& {R_0} \\
			{R_1} && {R_2} \\
			& {\{\abf\}}
			\arrow["{f_{1,0}}"', from=1-2, to=2-1]
			\arrow["{f_{2,0}}", from=1-2, to=2-3]
			\arrow[from=2-1, to=3-2]
			\arrow[from=2-3, to=3-2]
		\end{tikzcd}\]		
for some unital commutative rings $R_0$, $R_1$, and $R_2$ and some ring homomorphisms $f_{1,0}$, $f_{2,0}$.
\end{exmp}
	
	Recall that if $X$ is a topological space, then the set $\Open(X)$ of open subsets of $X$ is a lattice with the partial order defined by inclusion, having maximum $X$ and minimum $\emptyset$. We denote the category $\Pd\bigl(\Open(X)\bigr)$ by $\Pd(X)$ whenever there is no danger of ambiguity.
Similarly, we denote by $\Md(X)$ the full subcategory of $\Pd(X)$ consisting of common meadows $M$ such that $0\cdot M$ is isomorphic to $\Open(X)$.

	\begin{definition}\label{D:Operations}
		Let $X$ be a topological space, and $\fc$ a sheaf of rings in $X$. In the set  $M_\fc:=\bigsqcup_{U\in \Open(X)}\fc(U)$, \emph{addition} is defined by 
  $$s+t:=\rho_{U\cap V,U}(s)+\rho_{U\cap V,V}(t)$$ 
  
  and \emph{multiplication} is defined by $$s\cdot t:=\rho_{U\cap V,U}(s)\cdot \rho_{U\cap V,V}(t),$$ 	where $s\in \fc(U)$ and $t\in \fc(V)$ for $U$ and $V$ open subsets of $X$.
	\end{definition}

	\begin{exmp}\label{E:Cinf}
		Let $X=\r$ be the set of real numbers with the usual topology. Consider the sheaf $C^\infty$ of continuous real-valued functions and let $U$ be an open subset of $\r$. Then $C^\infty(U)=\{f:U\rightarrow \r  \mid f\text{ is continuous}\}$ is a ring with the sum and product defined pointwise. With the operations defined in Definition~\ref{D:Operations}, the set $M=\bigsqcup_{U\in \Open(X)}C^\infty(U)$ becomes a pre-meadow with $\abf$. Explicitly, given  $f\in C^\infty(U)$ and $g\in C^\infty(V)$, where $U$ and $V$ are open subsets of $\r$ we have
		\begin{enumerate}
			\item $f+g = f_{|U\cap V}+g_{|U\cap V}$
			\item $f\cdot g = f_{|U\cap V}\cdot g_{|U\cap V}$.
		\end{enumerate}
		
		Note that a function $f\in C^\infty(U)$ has an inverse in $C^\infty(U)$ if and only if $f$ is different from zero everywhere in $U$. Then its \emph{support} $$\supp(f):=\{x\in U\mid f(x)\neq 0\}$$ is an open subset of $\r$, and in fact the largest subset of $U$ where $f$ has an inverse. By Theorem \ref{T:DirectedLattice} we have that $M$ is a common meadow if and only if for all open subsets $U\subseteq X$ and all $f\in C^\infty(U)$ the set
		$$ J_{f}=\{V\in \Open(X)\mid f_{|V}\in C^\infty(V)^\times \}$$
		has a unique maximal element. We have seen that $\supp(f)\in J_{f}$. By the definition we have that $\supp(f)=\bigcup_{U\in J_f}U$ and so $\supp(f)$  is in fact the unique maximal element of $J_f$. Hence, $M$ is a common meadow.
	\end{exmp}

	The following lemma generalizes Example~\ref{E:Cinf}. As a consequence of the lemma, one obtains that morphisms of presheaves are the same as homomorphisms of pre-meadows with $\abf$. Also, and as illustrated in Examples~\ref{E:Projective}--\ref{E:Zariski} below, the functor $T$ defined in the lemma can be used to construct new kinds of meadows that, in particular, hint at a possible connection between meadows and algebraic geometry. This construction is extended to common meadows  in Theorem \ref{T:SheavesToCommon}.

	\begin{lemma}\label{L:Functor_T}
		
	Let $X$ be a topological space. Then we can define a functor
		\begin{align*}
			T:\PS&\rightarrow \Pd(X)\\
			\fc&\mapsto  M_{\fc}      
		\end{align*}
		by mapping each presheaf $\fc$ over $X$ to $M_{\fc}:=\bigsqcup_{U\in \Open(X)}\fc(U)$, the pre-meadow with $\abf$ endowed with the operations given in Definition \ref{D:Operations}, and mapping each morphism of presheaves ${\varphi : \fc \to \mathcal{L}}$ to the homomorphism of pre-meadows
		\begin{align*}
			\varphi':M_{\fc}&\rightarrow M_{\mathcal{L}}\\
			x&\mapsto  \varphi(U)(x),      
		\end{align*}
where $U$ is the unique open set such that ${x \in \fc(U)}$.
	\end{lemma}
	\begin{proof}
		 Recall that $\Open(X)$ is a lattice whose order is given by inclusion. We have that the set $\fc(U)$, where $U$ is an open subset  of $X$, is a unital commutative ring. Moreover, if $U\subseteq V$ are open subsets of $X$, or equivalently $V\leq U$ there exists a ring homomorphism: $\rho_{U,V}:\fc(V)\to\fc(U)$, then $(\fc(U),\rho_{U,V})_{U,V\in \Open(X)}$ is a directed lattice of rings. Hence $M$ is a pre-meadow with $\abf$, by Corollary \ref{C:DirectedLattice}.
		
		Let $\fc$ and $\mathcal{L}$ be two presheaves over $X$ and $\varphi:\fc\to \mathcal{L}$ a morphism of presheaves. Recall that for every $U\subseteq V$ open subsets of $X$ we have the following commutative diagram 
		\[\begin{tikzcd}
			{\mathcal{F}(V)} & {\mathcal{L}(V)} \\
			{\mathcal{F}(U)} & {\mathcal{L}(U)}
			\arrow["{\varphi(V)}", from=1-1, to=1-2]
			\arrow["{\rho_{U,V}^{\mathcal{F}}}", from=1-1, to=2-1]
			\arrow["{\rho_{U,V}^{\mathcal{L}}}", from=1-2, to=2-2]
			\arrow["{\varphi(U)}"', from=2-1, to=2-2]
		\end{tikzcd}\]
		We define the map
		$$ \varphi':T(\mathcal{F})\to T(\mathcal{L}) $$
		by $\varphi'(x):=\varphi(U)(x)$, where for each $x\in T(\mathcal{F})=\bigsqcup_{U\in \Open(X)}\fc(U)$, we are assuming that $x\in \fc(U)$, for some open subset $U$ of $X$. 
		
		Let us see that $\varphi'$ is an homomorphism of pre-meadows. Let $x\in \fc(U)$ and $y\in \fc(V)$. By definition, $x+y\in\fc(U\cap V)$. Then
  
		\begin{align*}
			\varphi'(x+y) &= \varphi(U\cap V)\bigl(\rho_{U\cap V,U}^{\mathcal{F}}(x)+\rho_{U\cap V,V}^{\mathcal{F}}(y)\bigr)\\
			&=\varphi(U\cap V)\bigl(\rho_{U\cap V,U}^{\mathcal{F}}(x)\bigr)+\varphi(U\cap V)\bigl(\rho_{U\cap V,V}^{\mathcal{F}}(y)\bigr)\\
			&=\rho_{U\cap V,U}^{\mathcal{L}}\circ\varphi(U)(x)+\rho_{U\cap V,V}^{\mathcal{L}}\circ\varphi(V)(y)\\
			&=\varphi'(x)+\varphi'(y).
		\end{align*}
  
		Similarly, and since $x\cdot y\in\fc(U\cap V)$, one shows that $\varphi'(x\cdot y)=\varphi'(x)\cdot \varphi'(y)$.
		
		Since $1_{T(\fc)}\in \fc(X)$, and  $\varphi(X)$ is a ring homorphism, we have that $\varphi(X)(1)$ is also the unity of the ring $\mathcal{L}(X)$. Hence $\varphi'$ is an homomorphism of pre-meadows.
It is straightforward to check that by making ${T(\varphi)=\varphi'}$, we have ${T(\Id_\fc)=\Id_{T(\fc)}}$ and ${T(\varphi \circ \psi)=T(\varphi) \circ T(\psi)}$ for any morphisms of presheaves $\varphi$ and $\psi$.
	\end{proof}

 As illustrated by the following example, finite common meadows can be obtained from the projective line, over any finite field. The example only considers the simplest case -- the one of the field $\z_2$ -- since for other finite fields the construction is similar, albeit more cumbersome to write and draw. 
 
	\begin{exmp}\label{E:Projective}
		Consider the field $\z_2$ and the sheaf of rings of regular functions on the projective line $\mathbb{P}^{1}$. Over $\z_2$ the projective line consists of three points, namely ${P_0=[0:1]}$, ${P_1=[1:1]}$, and ${P_\infty=[1:0]}$. We obtain the finite common meadow
		
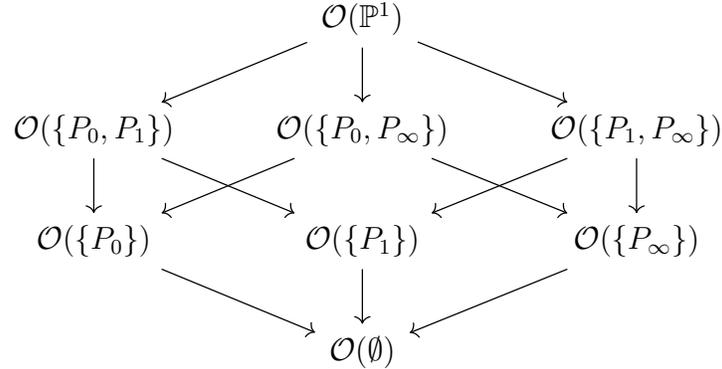
\begin{center}
			\begin{tikzcd}
				& {\mathcal{O}(\mathbb{P}^1)} \\
				{\mathcal{O}(\{P_0,P_1\})} & {\mathcal{O}(\{P_0,P_\infty\})} & {\mathcal{O}(\{P_1,P_\infty\})} \\
				{\mathcal{O}(\{P_0\})} & {\mathcal{O}(\{P_1\})} & {\mathcal{O}(\{P_\infty\})} \\
				& {\mathcal{O}(\emptyset)}
				\arrow[from=1-2, to=2-1]
				\arrow[from=1-2, to=2-2]
				\arrow[from=1-2, to=2-3]
				\arrow[from=2-1, to=3-1]
				\arrow[from=2-1, to=3-2]
				\arrow[from=2-2, to=3-1]
				\arrow[from=2-2, to=3-3]
				\arrow[from=2-3, to=3-2]
				\arrow[from=2-3, to=3-3]
				\arrow[from=3-1, to=4-2]
				\arrow[from=3-2, to=4-2]
				\arrow[from=3-3, to=4-2]
			\end{tikzcd}
		\end{center}
	\end{exmp}
 
	Let $\fc$ be a presheaf  over a topological space $X$. As shown in Example \ref{E:Notameadow} below,  the pre-meadow $M_\fc$ is not necessarily a common meadow. However, we show in Theorem \ref{T:SheavesToCommon} that if one starts with a sheaf instead of a presheaf, then the associated pre-meadow with $\abf$ is indeed a common meadow.
	
	\begin{exmp}\label{E:Notameadow}
		Let $X=\{a,b\}$ be a space with two points equipped with the discrete topology, so that every subset is open. Consider the presheaf $\fc$ on $X$ defined by $\fc(X)=\z$ and $\fc(\{a\})=\fc(\{b\})=\q$, and the restriction maps are given by the inclusion $i:\z\hookrightarrow \q$. 
		
		The presheaf $\q$ is not a sheaf since it does not satisfy the glueing property. In order to see this it is enough to consider $\frac{1}{3}\in \fc(\{a\})=\q $ and $\frac{1}{3}\in \fc(\{b\})=\q$. Indeed, there is no $m\in \fc(X)=\z$ that glues this two sections.
		
		The pre-meadow with $\abf$ associated with this presheaf is the one given by the following directed lattice:
		
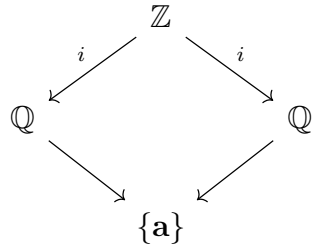
\[\begin{tikzcd}
			& \z \\
			\q && \q \\
			& {\{\abf\}}
			\arrow["i"', from=1-2, to=2-1]
			\arrow["i", from=1-2, to=2-3]
			\arrow[from=2-1, to=3-2]
			\arrow[from=2-3, to=3-2]
		\end{tikzcd}\]
		
		Given $m\in \z\setminus\{-1,0,1\}$ we have that $J_m$ has two distinct maximal elements which, by Theorem \ref{T:DirectedLattice}, shows that this pre-meadow with $\abf$ is not a common meadow. 
	\end{exmp}

The functor  $T$ from Lemma \ref{L:Functor_T} can be used to construct new kinds of common meadows where the directed lattice is infinite and related with regular functions over the projective line. Of course, similar constructions can also be made in higher dimensions, e.g.\ over the projective plane. 

	\begin{exmp}\label{E:Zariski}
		Consider the projective line $\mathbb{P}^{1}$ over a field $\mathbb{K}$ equipped with the Zariski topology. Since the closed sets are sets of zeroes of binary forms, they are either finite or the whole projective line. We may consider the sheaf of rings of regular functions on $\mathbb{P}^{1}$, i.e. functions that are defined locally as quotients of binary forms of the same degree. Applying the functor $T$ from Lemma \ref{L:Functor_T}, we obtain the following common meadow:
		
		\bigskip
		
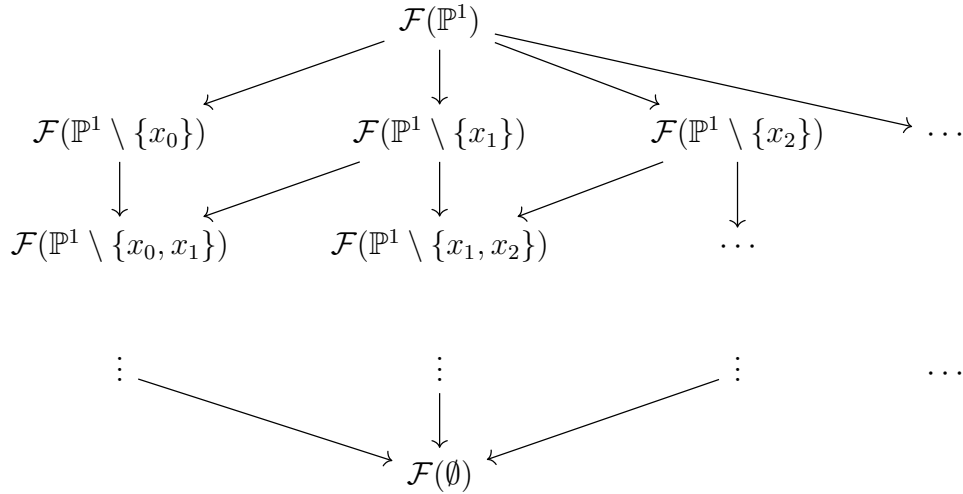
\begin{center}
			\begin{tikzcd}
				& {\mathcal{F}(\mathbb{P}^1)} \\
				{\mathcal{F}(\mathbb{P}^1\setminus\{x_0\})} & {\mathcal{F}(\mathbb{P}^1\setminus\{x_1\})} & {\mathcal{F}(\mathbb{P}^1\setminus\{x_2\})} & \cdots \\
				{\mathcal{F}(\mathbb{P}^1\setminus\{x_0,x_1\})} & {\mathcal{F}(\mathbb{P}^1\setminus\{x_1,x_2\})} & \cdots \\
				\vdots & \vdots & \vdots & \cdots \\
				& {\mathcal{F}(\emptyset)}
				\arrow[from=1-2, to=2-1]
				\arrow[from=1-2, to=2-2]
				\arrow[from=1-2, to=2-3]
				\arrow[from=1-2, to=2-4]
				\arrow[from=2-1, to=3-1]
				\arrow[from=2-2, to=3-1]
				\arrow[from=2-2, to=3-2]
				\arrow[from=2-3, to=3-2]
				\arrow[from=2-3, to=3-3]
				\arrow[from=4-1, to=5-2]
				\arrow[from=4-2, to=5-2]
				\arrow[from=4-3, to=5-2]
			\end{tikzcd}
		\end{center}
		Note that the vertical dots in the diagram are necessary because in case $\mathbb{K}$ is infinite, there is no minimal non-empty open set.
	\end{exmp}

	\begin{theorem}\label{T:SheavesToCommon}
		Let $X$ be a topological space, and $\fc$ a sheaf of rings on $X$. Then $M:=\bigsqcup_{U\in \Open(X)} \mathcal{F}(U)$ is a common meadow with the operations defined in Definition \ref{D:Operations}.
	\end{theorem}
	\begin{proof}
		By Lemma \ref{L:Functor_T} we know that $M$ is a pre-meadow with $\abf$. So we only need to check that, for every open set $W\subseteq X$ and  $f\in \fc(W)$ the set 
		$$
		J_f=\{ U\subseteq W \mid f_{|U}\in \fc(U)^\times\}
		$$
		has a unique maximal element.
		
		In order to do that it is enough to show that $V=\bigcup_{U\in J_f}U$ is an open set such that the restriction $f_{|V}$ is invertible. Since the elements of $J_f$ are open subsets we have that $V$ is also an open subset. From the definition of $J_f$ we have that for all $U\in J_f$ there exists $g_U\in\mathcal{F}(U)$ such that 
		
		\begin{align}
			(f_{|U}\cdot g_{U_1})_{|U_1\cap U_2}= f_{|U_1\cap U_2}\cdot (g_{U_1})_{|U_1\cap U_2} =1_{U_1\cap U_2}\\
			(f_{|U}\cdot g_{U_2})_{|U_1\cap U_2}= f_{|U_1\cap U_2}\cdot (g_{U_2})_{|U_1\cap U_2} =1_{U_1\cap U_2}.
		\end{align}
		By the uniqueness of the inverse, in the ring $\mathcal{F}(U_1\cap U_2)$ we have that $(g_{U_1})_{|U_1\cap U_2}=(g_{U_2})_{|U_1\cap U_2}$. Then, by the glueing property of the sheaf, there exists $g\in \mathcal{F}(V)$ such that $g_{|U}=g_U$ for all $U\in J_f$. It is straightforward to show that $g$ is the inverse of $f_{|V}$. Indeed, for all $U\in J_f$ it holds that
		$$
		(f_{|V}\cdot g)_{|U}= f_{|U}\cdot g_U=1_U.
		$$
		Since $J_f$ is by definition an open cover for $V$ we have that $f_{|V}\cdot g=1_V$ and so $V\in J_f$.
		
		Clearly, for any other maximal element $V'$ in $J_f$ we have $V'\subseteq \bigcup_{U\in J_f}U=V$, and since $V$ is in $J_f$ and $V'$ is maximal we also have $V\subseteq V'$. Hence $V$ is the unique maximal element of $J_f$.
	\end{proof}
	
	Given the nature of the functor $T$ it might seem simple to define an inverse functor. However, if one tries to do that, there is some ambiguity in choosing the presheaves, as illustrated in Example~\ref{E:ambiguity_sheaves}. This ambiguity may be overcome, and  one can indeed obtain an equivalence, by restricting the underlying topological space, as shown in Theorem~\ref{P:functorhomeo} below. 
	
	\begin{exmp}\label{E:ambiguity_sheaves}
		Consider the topological space $X=\{a,b\}$ with the discrete topology. The dual lattice of the open sets is the following
		
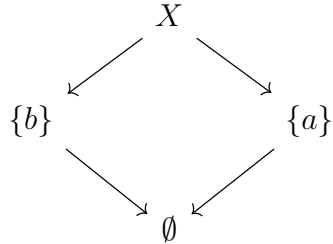
\[\begin{tikzcd}
			& X \\
			{\{b\}} && {\{a\}} \\
			& \emptyset
			\arrow[from=1-2, to=2-1]
			\arrow[from=1-2, to=2-3]
			\arrow[from=2-1, to=3-2]
			\arrow[from=2-3, to=3-2]
		\end{tikzcd}\]
		Then the pre-meadow with $\abf$ associated with the lattice
		
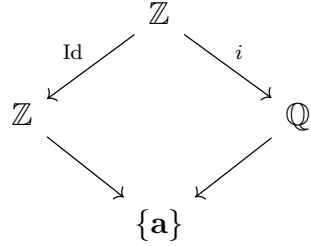
\[\begin{tikzcd}
			& \z \\
			\z && \q \\
			& {\{\abf\}}
			\arrow["\Id"', from=1-2, to=2-1]
			\arrow["i", from=1-2, to=2-3]
			\arrow[from=2-1, to=3-2]
			\arrow[from=2-3, to=3-2]
		\end{tikzcd}\]
		is clearly an element of $\Pd(X)$. To this pre-meadow with $\abf$ we can associate two different presheaves $\fc$ and $\fc'$ defined by:
		\begin{enumerate}
			\item $\fc(X)=\fc'(X)$
			\item $\fc(\{a\})=\z=\fc'(\{b\})$
			\item $\fc(\{b\})=\q=\fc'(\{a\})$
		\end{enumerate}
		
		Clearly these sheaves are not isomorphic, because otherwise we would have $\fc(\{a\})\simeq \fc'(\{a\})$.
	\end{exmp}

 \subsection{Relation with Sheafification}\label{S:Relation with Sheafification}
	Recall that a topological space $X$ is $T_D$ if all the singletons are locally closed, that is, if for all $x\in X$ the singleton $\{x\}$ is the intersection of an open subset of $X$ with a closed subset of $X$. Note that in particular Hausdorff spaces are $T_D$.  We recall the following result by W. J. Thron.
	
	\begin{theorem}\label{P:functorhomeo}\cite{thron1962lattice}
		Let $X$ and $Y$ be $T_D$ topological spaces. Then the following are equivalent:
		\begin{enumerate}
			\item There exists an homeomorphism $f:X\to Y$.
			\item There exists a lattice homomorphism $f':\Open(X)\to \Open(Y)$.
		\end{enumerate}
		Moreover, the maps are compatible, meaning that given $U$ an open subset of $X$ we have $f'(U)= f(U)$. 
	\end{theorem}
	
	Our goal now is, starting with a pre-meadow with $\abf$, to construct an associated common meadow making use of the sheafification of a certain presheaf, defined in a $T_D$ topological space.
	
	\begin{lemma}\label{L:Functor_T'}
		Let $X$ be a topological space, and take $P\in \Pd(X)$. Let $\Phi:0\cdot P\to \Open(X)$ be an isomorphism of lattices. Then, there exists a presheaf of rings $\fc_{P,\Phi}$ on $X$ defined by $\fc_{P,\Phi}(U)=P_{\Phi\inv (U)}$.
	\end{lemma}
	\begin{proof}
		This is an immediate consequence of the definition of presheaf.
	\end{proof}
	
	\begin{lemma}\label{L:SheafIsom}
		Let $X$ be a $T_D$ topological space and let $\fc$ and $\fc'$ be two presheaves in $X$. Then $T(\fc)\simeq T(\fc')$ if and only if there exists an homeomorphism $f:X\to X$ such that $\fc\simeq f_*(\fc')$, where $f_*$ is the direct image functor.
	\end{lemma}
	\begin{proof}
		By construction we have that $T(\fc)=\bigsqcup_{U\in \Open(X)}\fc(U)$ and $T(\fc')=\bigsqcup_{U\in \Open(X)}\fc'(U)$.
		
		Assume first that $\Phi:T(\fc)\to T(\fc')$ is an  isomorphism of meadows. Then, we also have  an isomorphism of lattices, defined by
		$$\Phi_L:\Open(X)=0\cdot T(\fc)\to 0\cdot  T(\fc)= \Open(X).$$
		By Theorem \ref{P:functorhomeo}, the isomorphism $\Phi_L$ corresponds to an homeomorphism $\Phi':X\to X$ such that,  given $V\subseteq U$ open subsets of $X$, the following diagram commutes
		
		\[
		\begin{tikzcd}
			{\fc(U)=T(\fc)_U} && {T(\fc')_{\Phi_L(U)}=\fc'(\Phi'(U))} \\
			\\
			{\fc(V)=T(\fc)_V} && {T(\fc')_{\Phi_L(V)}=\fc'(\Phi'(V))}
			\arrow["{\Phi_U}", from=1-1, to=1-3]
			\arrow["{f_{V,U}^{T(\fc)}}"', from=1-1, to=3-1]
			\arrow["{f_{V,U}^{T(\fc')}}", from=1-3, to=3-3]
			\arrow["{\Phi_V}"', from=3-1, to=3-3]
		\end{tikzcd}
		\]
		
		In the diagram above $\Phi_U$ is a ring isomorphism. Since $\Phi'$ is an homeomorphism we have that the commutativity of the diagram is exactly the condition for the sheaves $\fc$ and $(\Phi')_*(\fc')$ to be isomorphic. 
		
		Assume now that there exists an homeomorphism $f:X\to X$ such that there exists $\varphi:\fc\to f_*(\fc')$ an isomorphism of sheaves. By definition, given $U\subseteq V$ open subsets of $X$ we have the following diagram

		\[\begin{tikzcd}
	{\fc(U)} && {(f_*\fc')(U)=\fc'(f\inv(U))} \\
	\\
	{\fc(V)} && {(f_*\fc')(V)=\fc'(f\inv(V))}
	\arrow[from=1-1, to=1-3]
	\arrow["{{\varphi(U)}}", from=1-1, to=1-3]
	\arrow["{{\rho_{V,U}^{\fc}}}", from=1-1, to=3-1]
	\arrow["{{\rho_{V,U}^{f_*\fc'}=\rho_{f\inv(V),f\inv(U)}^{\fc'}}}", from=1-3, to=3-3]
	\arrow["{{\varphi(V)}}", from=3-1, to=3-3]
\end{tikzcd}\]
		
		By the same argument as before we have that $\overline{\varphi}:T(\fc)\to T(\fc')$ defines an isomorphism of common meadows, where given an $s\in \fc(U)\subseteq T(\fc)$ we define $\overline{\varphi}(s)=\varphi(U)(s)\in \fc'(f\inv(U))\subseteq T(\fc')$.
	\end{proof}
	
	\begin{theorem}\label{T:Common_Meadow_Sheafification}
		Let $X$ be a $T_D$ topological space, $P\in \Pd(X)$, and $\Phi:0\cdot P\to \Open(X)$ be an isomorphism of lattices. Then there is an associated common meadow defined by $M:=T(\fc_{P,\Phi}^*)$, where  $\fc_{P,\Phi}^*$ is the sheafification of the presheaf $\fc_{P,\Phi}$.
		
		Moreover, $M$ is independent of the choice of the isomorphism $\Phi$, in the sense that if $\Phi':0\cdot P\to \Open(X)$ is another isomorphism of lattices then $T(\fc_{P,\Phi}^*)$ is isomorphic to $T(\fc_{P,\Phi'}^*)$. 
	\end{theorem}
	\begin{proof}
		Given $P$ a pre-meadow with $\abf$, we may use Lemma \ref{L:Functor_T'} to obtain a sheaf $\fc_{P,\Phi}^*$. Then  $T(\fc_{P,\Phi}^*)$ is a common meadow by Theorem \ref{T:SheavesToCommon}.
		
		Let us see that this is in fact independent of the choice of $\Phi$. Note that $T(\fc_{P,\Phi})$ and $T(\fc_{P,\Phi'})$ are isomorphic by construction. Then, by Lemma \ref{L:SheafIsom}, there exists an homeomorphism $f:X\to X$ such that $\fc_{P,\Phi}\simeq f_* (\fc_{P,\Phi'})$. Since, in this case, the direct image functor commutes with the sheafification (Remark~\ref{R:directimagecommutes}) we have that $\fc_{P,\Phi}^*\simeq f_* ((\fc_{P,\Phi'})^*)$. Then, Lemma \ref{L:SheafIsom} implies that $T(\fc_{P,\Phi}^*)$ is isomorphic to $T(\fc_{P,\Phi'}^*)$, as we wanted.
	\end{proof}
	
	Note that if $M$ is a common meadow and $\Phi:0\cdot M\to \Open(X)$ is an isomorphism of lattices, the common meadow $T(\fc_{M,\Phi}^*)$ is not necessarily isomorphic to $M$ as illustrated by the following example.
	
	\begin{exmp}
		Consider $X=\{a,b\}$ equipped with the discrete topology. Then the pre-meadows $P$ in $\Pd(X)$ have a lattice isomorphic to the following lattice
		
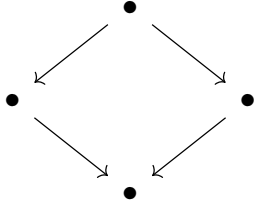
\[\begin{tikzcd}
			& \bullet \\
			\bullet && \bullet \\
			& \bullet
			\arrow[from=1-2, to=2-1]
			\arrow[from=1-2, to=2-3]
			\arrow[from=2-1, to=3-2]
			\arrow[from=2-3, to=3-2]
		\end{tikzcd}\]
		Consider the common meadow $M$ defined by the directed lattice
		
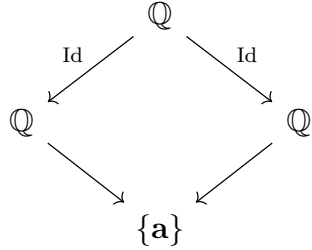
\[\begin{tikzcd}
			& \q \\
			\q && \q \\
			&\{\abf\}
			\arrow["\Id"', from=1-2, to=2-1]
			\arrow["\Id", from=1-2, to=2-3]
			\arrow[from=2-1, to=3-2]
			\arrow[from=2-3, to=3-2]
		\end{tikzcd}\]
		Then, the associated presheaf is the presheaf $\fc_M$ such that for any non-empty open subset $U\subseteq X$ we have $\fc_M(U)=\q$. And the sheafification $\fc_M^*$ of the presheaf $\fc_P$ is $\fc_M^*(\{a\})=\fc_M^*(\{b\})=\q$, and $\fc_M^*(X)=\q\oplus \q$. Then the common meadow $T(\fc_M^*)$ is defined by the following lattice
		
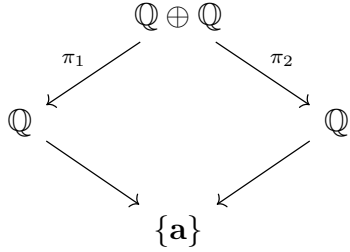
\[\begin{tikzcd}
			& \q\oplus\q \\
			\q && \q \\
			& \{\abf\}
			\arrow["{\pi_1}"', from=1-2, to=2-1]
			\arrow["{\pi_2}", from=1-2, to=2-3]
			\arrow[from=2-1, to=3-2]
			\arrow[from=2-3, to=3-2]
		\end{tikzcd}\]
		which is not isomorphic to $M$, since $M_0=\q$ and $(T(\fc_M^*))_0=\q\oplus\q$.
	\end{exmp}

 \subsection{An isomorphism}\label{S:Isomorphism}
 Theorem~\ref{T:Common_Meadow_Sheafification} shows that the process of going from a pre-meadow with $\abf$ to a common meadow can be identified with the sheafification of a presheaf. However, as illustrated by Example~\ref{E:Not_Iso}, if one already starts with a common meadow, the process may give a common meadow not isomorphic to the original one.

 \begin{exmp}\label{E:Not_Iso}
    Consider the constant pre-sheaf $\underline{\z}$ on the set of real numbers $\r$, where to each non-empty subset we associate the ring of integers $\z$, and to the empty set the trivial ring. As seen in Example~\ref{E:Not_a_Sheaf}, the presheaf $\underline{\z}$ is not a sheaf. However, the associated pre-meadow $T(\underline{\z})$ is a common meadow, since for every $0\cdot z\in 0\cdot T(\underline{\z})\setminus \{\abf\} $ we have that $T(\underline{\z})_{0\cdot z}=\z$ (see \cite[Corollary~2.10]{Dias_Dinis(Art)}). However, the common meadow $T(\underline{\z}^*)$ is not isomorphic to $T(\underline{\z})$ since $\underline{\z}^*$ is not constant.
    \end{exmp}

	The conditions for a presheaf to be a sheaf can nevertheless be adapted to give conditions to establish an isomorphism between a class of common meadows $M$ and $T(\fc_M^*)$. First recall that given an isomorphism of lattices $\Phi:0\cdot M\to \Open(X)$ there exists a join operation $\vee$ in $0\cdot M$ induced by the union of open subsets.
	
	\begin{theorem}
		\label{P:functorequi}
		
		Let $X$ be a topological space and $M$ be a common meadow such that  ${\Phi:0\cdot M\to \Open(X)}$ is an isomorphism of latices. Then $T(\fc_{M,\Phi}^*)$ is isomorphic to $M$ if and only if for all $0\cdot z,0\cdot w\in 0\cdot M$ the following two conditions hold:
		\begin{enumerate}
			\item \label{P:functorequi1}If $x,y\in M_{(0\cdot z)\vee (0\cdot w)}$ are such that 
			$\begin{cases}
			{x+0\cdot   z =y+0\cdot   z}\\ 
			{x+0\cdot   w =y+0\cdot  w},
			\end{cases}$
			then $x=y$.
			
			\item \label{P:functorequi2} 
			If $\{x_i\}_{i\in I}$ is a collection of elements of $M$ such that for any ${i,j\in I}$ it holds that
			\[
			x_i+0\cdot x_i\cdot x_j=x_j+0\cdot x_i\cdot x_j,
			\] 
			then there exists some ${x\in M}$ such that ${x_i=x+0\cdot x_i}$, for all ${i\in I}$.
		\end{enumerate}
	\end{theorem}
	
	\begin{proof}
            We will first prove that $\fc_{M,\Phi}$ is a sheaf if and only if $M$ satisfies conditions \ref{P:functorequi1} and \ref{P:functorequi2}.

            From the definition of sheaf we have that $\fc_{M,\Phi}$ is a sheaf if and only if
            \begin{enumerate}
			\item[$(a)$]
For any $s,t\in \fc_{M,\Phi}(U)$, if $\rho_{U_i,U}(s)=\rho_{U_i,U}(t)$ for all $i\in I$ then $s=t$.
			\item[$(b)$] 
For any collection $\{s_i\}_{i\in I}$ satisfying ${s_i\in \fc_{M,\Phi}(U_i)}$ for each ${i \in I}$, if we have $\rho_{U_i\cap U_j,U_i}(s_i)=\rho_{U_i\cap U_j,U_j}(s_j)$ for all ${i,j\in I}$ then there exists $s\in U$ such that $\rho_{U_i,U}(s)=s_i$, for all ${i \in I}$.
		\end{enumerate}
        where $\rho_{U_i,U}$ are the restriction maps of the sheaf  $\fc_{M,\Phi}$. Since the restriction map $\fc_{M,\Phi}(U)\to \fc_{M,\Phi}(U_i)$ is given by $\rho_{U,U_i}(x)=x+\Phi\inv(U_i)$ we have that condition $(a)$ and $(b)$ are the same as conditions \ref{P:functorequi1} and \ref{P:functorequi2}, respectively.
        
        If  $\fc_{M,\Phi}$ is a sheaf, then clearly $T(\fc_{M,\Phi})=M$, as we wanted.

        Suppose now that there is an isomorphism between $T(\fc_{M,\Phi}^*)$ and $M$. Then $\fc_{M,\Phi}^*$ is isomorphic to $\fc_{M,\Phi}$. Since $\fc_{M,\Phi}^*$ is a sheaf we have that $\fc_{M,\Phi}$ must also be a sheaf, and so conditions  \ref{P:functorequi1} and \ref{P:functorequi2} are satisfied.
	\end{proof}

	
	\section{Conclusion}\label{S:Conclusion}

	This paper established a connection between meadows and algebraic topology by providing a construction that given a sheaf of rings $\mathcal{F}$ on a topological space $X$ produces a common meadow as a disjoint union of elements of the form $\mathcal{F}(U)$ indexed over the open subsets of $X$ (Theorem~\ref{T:SheavesToCommon}). In the other direction, we gave a correspondence between the process of sheafification and the process of going from a pre-meadow with $\abf$ to a common meadow (Theorem~\ref{T:Common_Meadow_Sheafification}). Of course, now that such a connection is made one may ask how can this connection be explored. For example,
	what happens when we use other well-known topologies? Are there some interesting ``extreme'' cases?

 The paper \cite{Dias_Dinis(24)} introduced an algorithm to count the number of finite  pre-meadows with $\abf$. As mentioned in that same paper, the ultimate goal is to obtain a function enumerating all finite common meadows. It is indeed not straightforward to expand an enumeration of pre-meadows with $\abf$ to an enumeration of common meadows, but perhaps the connection with sheaves presented here can be a step towards an answer to this question. More generally,  we hope that the connections presented in this paper can enable one to translate problems from algebraic topology to the setting of meadows (or the other way around) where they might become more tame. So far we have only scratched the surface and used tools from algebraic topology to construct new kinds of meadows, as can be seen for instance in Examples \ref{E:Projective} and \ref{E:Zariski}. These examples also hint at a possible deeper connection with algebraic geometry. 
	
\subsection*{Acknowledgments}

All authors acknowledge the support of FCT - Funda\c{c}\~ao para a Ci\^{e}ncia e Tecnologia under the project:\\ 10.54499/UIDB/04674/2020, and the research center CIMA -- Centro de Investigação em Matemática e Aplicações. 

The second author also acknowledges the support of CMAFcIO -- Centro de Matem\'{a}tica, Aplica\c{c}\~{o}es Fundamentais e Investiga\c{c}\~{a}o Operacional under the project UIDP/04561/2020.


\end{document}